\documentclass[12pt]{amsart}
\usepackage{amsmath,amssymb}
\usepackage{latexsym}
\usepackage{amsthm}
\title[Pseudoreflection Action]{Groups of Finite Morley Rank with a Pseudoreflection Action}
\author{Ay\c{s}e Berkman}
\address{Mathematics Department, Mimar Sinan University, Silah\c{s}\"or Cad. 71, Bomonti \c{S}i\c{s}li 34380, \.Istanbul, Turkey.}
\email{ayse.berkman@msgsu.edu.tr}
\author{Alexandre Borovik}
\address{School of Mathematics,  University of Manchester, Oxford Road, Manchester M13 9PL}
\email{borovik@manchester.ac.uk}

\begin{document}
\newtheorem{problem}{Problem}
\newtheorem{lemma}{Lemma}[section]
\newtheorem{theorem}[lemma]{Theorem}
\newtheorem{corollary}[lemma]{Corollary}
\newtheorem{proposition}[lemma]{Proposition}
\newtheorem{remark}[lemma]{Remark}
\newtheorem{fact}[lemma]{Fact}
\newtheorem{question}[lemma]{Question}
\theoremstyle{definition}
\newtheorem*{definition}{Definition}
\newcommand{\acf}{algebraically closed field }
\newcommand{\acfd}{algebraically closed field}
\newcommand{\acfs}{algebraically closed fields}
\newcommand{\fmr}{finite Morley rank }
\newcommand{\rk}{{\rm rk}}
\newcommand{\psrk}{{\rm psrk}}
\newcommand{\fmrd}{finite Morley rank}
\newcommand{\bi}{\begin{itemize}}
\newcommand{\ei}{\end{itemize}}
\maketitle

\begin{abstract}
In this work, we give two characterisations of the general linear group as a group $G$ of \fmr acting on an abelian connected group $V$ of \fmr definably, faithfully and irreducibly. To be more precise, we prove that if the pseudoreflection rank of $G$ is equal to the Morley rank of $V$, then $V$ has a vector space structure over an \acfd, $G\cong \operatorname{GL}(V)$ and the action is the natural action. The same result holds also under the assumption of Pr\"ufer 2-rank of $G$ being equal to the Morley rank of $V$.
\end{abstract}
\section{Introduction}

Groups of finite Morley rank are abstract groups
equipped with a notion of dimension which assigns to every
definable set $X$ a dimension, called {\em Morley rank} and
denoted ${\rm rk}(X)$, satisfying well known and fairly
rudimentary axioms given for example in \cite{abc,bn,poizat}.
Examples are furnished by algebraic groups over algebraically
closed fields, with ${\rm rk}(X)$ equal to the dimension of the
Zariski closure of $X$.

Crucially, groups of finite Morley rank naturally arise in model theory. By work of Boris Zilber, any uncountably categorical structure is controlled by certain definable groups of permutations (which have finite Morley rank, by definability). This observation leads to the concept of a \emph{binding group}, introduced by Zilber and developed in other contexts by Poizat and Hrushovski. Binding groups play in model a role akin to that  of Galois
groups; Lie groups of the Picard-Vessiot theory of linear differential equations can be viewed as a special case, as it had been explained by Poizat \cite{Poizat83}. The following theorem gives an example of appearance of a binding group:

\begin{fact}[\cite{poizat}]
Let $T$ be an $\omega$-stable theory, $M \models T$ a prime model  over $\emptyset$, and let $P$, $Q$ be
$\emptyset$-definable sets, with $P$ being $Q$-internal {\rm (}that is, $P \subset  {\rm dcl}(Q \cup F)$ for some finite $F${\rm )}. Then the group of automorphisms of $M$ which fix $Q$ pointwise induces a
definable group of automorphisms of $P$, called the \emph{binding group} of $P$ over
$Q$.
\end{fact}

Analysing relations between model theory and group theory, one sees, first of all, the profound impact of the Classiffication of Finite Sumple Groups on development of model theory, achieved mostly via finite combinatorics and the theory of permutation groups.

Even without the complete classification of groups of finite Morley rank,
the results and methods already developed and systematised in works like \cite{abc} are powerful enough to
start a systematic structural theory of binding groups, first in the context of finite Morley rank and then in more general context of stable and simple theories.  This will require a good understanding of permutation groups of finite Morley rank which makes a natural first stage of the project.

In \cite{borche}, Borovik and Cherlin stated several problems regarding permutation groups of \fmrd. Among other things, they asked about generic multiply transitive actions of groups of \fmrd, and as a first step towards a full description, they asked the following question:

Assume that $G$ is a connected group of \fmr acting on a
connected abelian group $V$ of \fmr definably, faithfully and
irreducibly. Also $G$ contains a pseudoreflection subgroup $R$,
that is a connected definable abelian subgroup such that $V=[V,R]\oplus
C_V(R)$, and $R$ acts transitively on the nonzero elements of
$[V,R]$. Then is it true that $V$ is a vector space and $G\cong \operatorname{GL}(V)$?

This paper gives an affirmative answer to this question under the extra assumptions that $\rk[V,R]=1$ and $\rk(V)=\psrk(G)$ where $\psrk(G)$ is the \emph{pseudoreflection rank} of $G$, that is, the maximal number of pairwise commuting pseudoreflection subgroups in $G$. The reader may note that these rank assumptions do not put any restriction for the ``real life'' case, where the identification of the group is possible. The remaining case is a non-existence proof, so far it seems as if it will be a long technical discussion about non-existing objects.

In this paper, we  prove the following.

\begin{theorem}\label{main} Let $G$ be a connected group of \fmr acting on a
connected abelian group $V$ of \fmr definably, faithfully and
irreducibly. If\/ $G$ contains a pseudoreflection subgroup $R$ such that\/ $\rk[V,R]=1$, and\/  $\rk(V)=\psrk(G)$, then $V$ is a vector space over an \acfd, $G\cong \operatorname{GL}(V)$ and acts on $V$ as on its natural module. \end{theorem}

The following corollary easily follows from the above theorem.

\begin{corollary} Let $F$ be an \acfd, and $G$ a connected group of finite Morley rank acting definably and faithfully on $F^n$ by automorhisms of $F^n$ as an  additive group {\rm (}but not necessarily preserving the structure of the $F$-vector space{\rm )}. If\/ $GL_n(F)$ lies in $G$, then $G=GL_n(F)$.
\end{corollary}

We also deduce the following theorem whose proof will be given in Section~\ref{sec:large-torus}.

\begin{theorem}\label{large-torus}
Let $G$ be a connected group of \fmr acting on a
connected abelian group $V$ of \fmr definably, faithfully and
irreducibly. If\/ ${\rm pr_2}(G) = \rk(V)$  then $V$ is a vector space over an \acfd, $G\cong \operatorname{GL}(V)$ and acts on $V$ as on its natural module.
\end{theorem}

The outline of our proof is as follows. First, we use centraliser functors to prove that for every non-central involution $i$ in $G$, $C_G^\circ(i)$ is isomorphic to a direct sum of copies of $\operatorname{GL}_{n_i}$, where $\sum_{i} n_i=\rk(V)$ (Section~\ref{reduct}); this allows us to apply the Generic Identification Theorem \cite{git,git2} and complete the proof of Theorem~\ref{main} (Section~\ref{sec;git}).

\section{Background Material}

From now on, all groups are of \fmrd. We will be using \cite{bn} or \cite{abc} for general reference in this work.

Since in a group $G$ of \fmrd, definable subgroups satisfy the descending chain condition \cite[Theorem 5.2]{bn}, $G$ contains a minimal definable subgroup of finite index, which is called the {\em connected component} of $G$, and is denoted by $G^\circ$.
If $G=G^\circ$, then we say $G$ is {\em connected}.

If $H$ is a (not necessarily definable) subgroup of $G$, then we define the \emph{connected component} of $H$ as  $H^\circ=d(H)^\circ$, where $d(H)$ is the smallest definable subgroup in $G$ containing $H$, it is called the {\em definable hull} of $H$.

A divisible abelian group is called a {\em torus}. For a prime $p$, a $p$-group which is also a torus is called {\em $p$-torus}. Hence, a $p$-torus is a direct product of copies of the quasi-cyclic group $C_{p^\infty}$.

\begin{itemize}
\item A subgroup $H$ of $G$ is {\em unipotent} if it is definable,
  connected, nilpotent, and of bounded exponent.
\item $G$ is a {\em $K$-group} if every infinite connected simple definable
section of $G$ is a Chevalley group over an algebraically closed field.

\item A {\em Sylow$^\circ$ $2$-subgroup} of $G$ is the connected component
  of a maximal $2$-subgroup.
\item $G$ is of {\em even type} if its Sylow$^\circ$ $2$-subgroup
is nontrivial and unipotent.

\end{itemize}

The following result of Zilber is very useful.

\begin{fact} [Zilber, Theorem 9.1 \cite{bn}] If an infinite abelian (solvable) group $A$ acts on an infinite abelian group $V$ definably, faithfully and minimally, then $A$ is isomorphic to a subgroup of $K^*$, $V\cong K^+$ for some \acf $K$, and $A$ acts on $V$ by multiplication. \label{zilberabelianaction}
\end{fact}

\begin{fact} [Poizat]  If $F$ is an \acf of positive characteristic, then any infinite definable simple subgroup of $GL_n(F)$ is an algebraic group over $F$.\label{poizat}
\end{fact}

\medskip
\noindent
{\bf Remark.} Poizat's original proof of the above fact uses the classification of finite simple groups, see \cite[Proposition II.4.4]{abc} for a sketch of the proof. However, Borovik supplied another proof in \cite{avbrussian} without using the classification of finite simple groups.

\begin{fact} {\rm \cite{bc}} Let $G$ be a connected group of \fmrd, and $p$ a prime number. If $a\in G$ is a $p$-element and $C_G^\circ(a)$ has no nontrivial $p$-unipotent subgroups, then $a$ lies in a $p$-torus. \label{burche}
\end{fact}

A torus $T$ is called {\em good} if every definable subgroup $A$ in $T$ is the definable hull of the torsion subgroup of $A$.

Note that if $F$ is an infinite field of finite Morley rank and
of positive characteristic, then $F^*$ is a good torus.

\begin{fact} In a group of \fmr $G$, maximal good tori in $G$ are conjugate in $G$. %Maximal decent tori are conjugate.
If $T$ is a good torus in $G$, then $N^\circ_G(T)=C^\circ_G(T)$. Moroever, if $G$ is connected, then $C_G(T)$ is connected.\label{cherlintorus}\label{altbur}
\end{fact}

\proof The first statements appeared in \cite{cherlin}, and the last statement is from \cite{altbur}. See also \cite[Chapter IV]{abc}.

\section{Special Cases of Theorem~\ref{main}}

First, we shall treat the $\rk(V)=1$ case which follows from the following lemma.

\begin{lemma} If $G$ is a connected group acting faithfully on a strongly minimal group $V$, then there exists an \acf $F$ such that $V\cong F^+$, $G\leqslant F^*$ and the action is the usual multiplication. \label{strmin}
\end{lemma}

\proof First let us prove the following.

\medskip
\noindent
{\em Claim 1.} If $G$ is a connected group acting non-trivially on a strongly minimal set $V$, then $C_V(G)$ is finite and $G$ acts transitively on the cofinite subset $V\setminus C_V(G)$.

\medskip
\noindent
{\em Proof of Claim 1.} Let $v\in V$, then its orbit $orb(v)$ under the action of $G$ is a definable subset of $V$, hence it is either finite or cofinite. If $orb(v)$ is finite, say it has $n$ elements, then the stabiliser $stab(v)$ of $v$ under the action of $G$ is a definable subgroup of $G$ of index $n$, hence $n=1$. Note that we have the following now: $orb(v)$ is finite iff $orb(v)=\{v\}$ iff $v\in C_V(G)$. Again $C_V(G)$ is a definable subset of $V$, hence it is finite or cofinite.
If $C_V(G)$ is cofinite, then the connected group $G$ acts on the finite set $V\setminus C_V(G)$, hence it acts trivially on it. Thus, $V=C_V(G)$, but this contradicts with the non-triviality of the action. Thus, $C_V(G)$ is finite. Therefore, $orb(v)$ is cofinite and $V=orb(v)\sqcup C_V(G)$ for every $v\in V\setminus C_V(G)$.

\medskip
\noindent
{\em Claim 2.} If a connected group $G$ acts on a strongly minimal group $V$ definably and faithfully, then $G$ has at most one involution.

\medskip
\noindent
{\em Proof of Claim 2.} Let $i$ be an involution in $G$. Since $C_V(i)$ is a definable subgroup of $V$, $C_V(i)$ is either finite or $C_V(i)=V$. The latter contradicts the faithfullness of the action. Hence, $C_V(i)$ is finite and by \cite[I.10.3]{abc}, $i$ acts on $V$ by inversion. Now the uniqueness of the involution (if exists) follows from the faithfullness of the action.

\medskip
\noindent
{\em Proof of the theorem.} Now, by Claim 1, there exists a cofinite subset $W\subseteq V$ on which $G$ acts transitively. Note that $W$ is also strongly minimal and the action of $G$ on $W$ is faithful. Hence we can apply the classification of transitive group actions on strongly minimal sets, done by Hrushovski \cite[Theorem 11.xx]{bn}, and conclude that $G$ is either solvable or $G\cong \operatorname{PSL}_2(F)$ over some \acf $F$.

If $G$ is solvable, then $G'$ acts trivially on $V$, hence $G$ is abelian. Thus we get the desired result by Fact~\ref{zilberabelianaction}. Since $\operatorname{PSL}_2(F)$ has more than one involution, $G\cong \operatorname{PSL}_2(F)$ is not possible by Claim~2. \hfill $\Box$

The case $\rk(V)=2$ is contained in a more general result by Deloro.

\begin{fact}[Deloro {\rm \cite{deloro}}] \label{deloro} Let $G$ be a connected non-solvable group of \fmr acting definably and faithfully on a connected abelian group $V$ of Morley rank 2. Then there exists an \acf $K$ such that $V\cong K^2$ and $G\cong \operatorname{SL}_2(K)$ or $\operatorname{GL}_2(K)$ in their natural action on $K^2$.\end{fact}

Two important special cases of Theorem~\ref{main} follow from the following two facts.

\begin{fact} [Zilber, Theorem 9.5 \cite{bn}] If a connected group $G$ with an infinite abelian normal subgroup acts on an abelian connected group $V$ faithfully and minimally, then $G$ lies in $GL(V)$.\label{zilberlin}
\end{fact}

\begin{fact} {\rm \cite{abc}} Under the assumptions of Theorem~{\rm\ref{main}}, if $V\cong F^n$, for some \acf $F$, and $G\leqslant \operatorname{GL}(V)$, then $G\cong\operatorname{GL}(V)$.\label{abclin}
\end{fact}

The proof of this fact is known, for sake of completeness, we extract it from the proof of Theorem III.1.5 in \cite{abc}.

\medskip

\proof By the above quoted facts, we may assume $n\ge 3$. Let $H$ stand for the subgroup in $G$ generated by all the pseudoreflection subgroups. Then $H$ is a connected definable algebraic group (since pseudoreflection subgroups are Zarsiki closed) %at most one eigenvalue not 1, simult. diagbl.
and is a normal subgroup in $G$.

Assume that $H$ acts reducibly on $V$. Then by the Clifford Theorem, $V=\bigoplus V_i$. Each reflection subgroup in $H$ acts on exactly one $V_i$ non-trivially. Hence, we can write $H=H_1\times\cdots\times H_k$ where $H_i$ is generated by those pseudoreflection subgroups that act on $V_i$ non-trivially. Now note that $V_1,\ldots, V_k$ are the only irreducible $H$-modules in $V$. Thus, the connected group $G$ acts transitively on them by Clifford and the previous statement. Hence there is only one $H$-module, which is $V$.

Therefore, $H$ acts irreducibly on $V$ and thus $H$ is reductive. (Recall that unipotent radicals fix a nontrivial subgroup, hence an irreducible group with pseudoreflections is reductive.) By the structure theorem for reductive algebraic groups, $H=TM$ an almost direct product of a torus $T$ and a product $M$ of quasisimple subgroups. Since the action is irreducible, $M\neq 1$. Note that $T$ consists of scalars by the Schur Lemma, that is, $T$ is a one-dimensional torus.

Now note that $H$ contains root subgroups. Thus, $K=\langle R,R^g\rangle$ acts faithfully on $W=[V,R]+[V,R^g]$, and hence $K\leqslant GL_2$. By computation, $K$, and hence $H$, contains root subgroups. It is easy to see that the subgroup in $H$ which is generated by root subgroups is $M$. Since $H=T\times M$ and $H$ acts irreducibly on $V$, $M$ also acts irreducibly on $V$. Thus, we are in the setting of McLaughlin's theorem.

\begin{fact} [McLaughlin \cite{mcl}] Let $F$ be a field with at least $3$ elements, $V$ a vector space of finite dimension at least $2$ over $F$, $1\neq M\leqslant \operatorname{SL}(V)$ generated by subgroups of root type acting irreducibly on $V$ and $R_u(M)=1$. Then, $M\cong \operatorname{SL}(V)$ or $\operatorname{Sp}(V)$.
\end{fact}

Therefore, $M\cong \operatorname{SL}(V)$ or $M\cong \operatorname{Sp}(V)$ in $H$. However, $\operatorname{Sp}(V)$ preserves a skew-symmetric form, but $T$ cannot preserve it. %More precisely, take $T$ and the maximal pseudoreflection subgroup in $M$ and .... obtain $diag(a,b,\ldots,b)$ where $ab^{n-1}=1$. This lies in $M$ but does not preserve the form. CHECK!!!
Hence, $M\cong \operatorname{SL}_n$ and $H=M\times T\cong \operatorname{GL}_n$. Now we have $\operatorname{GL}_n\cong H\leqslant G\leqslant \operatorname{GL}_n$, hence $G=\operatorname{GL}_n$.   \hfill $\Box$

\subsection{Charactersitic 0 and 2}

Other two known special cases of Theorem~\ref{main} are related to the structure of $V$.

From now on, $V$ is a connected abelian group, $G$ is connected and acts definably, faithfully and irreducibly on $V$ with a pseudoreflection subgroup $R$ such that $\rk[V,R]=1$, and $\rk(V)=n$.

First, note that an immediate corollary of Zilber's Theorem (Fact~\ref{zilberabelianaction}) in our setting is that $R\cong F^*$ and $[V,R]\cong F^+$ for some \acf $F$.

Note that by irreducibility $$V=\left\langle \bigcup_{g\in G}\left[V,R^g\right]\right\rangle =\sum_{i=1}^n\left[V,R^{g_i}\right]\cong\sum_n F^+.$$
The second equality follows from Zilber's Indecomposibility  Theorem \cite[5.28]{bn}.
(Moreover, when characteristic is 0, $V=\bigoplus_n F^+$, since $F^+$ has no nontrivial proper definable subgroups.) Therefore, $V$ is either a torsion-free divisible abelian group or an elementary abelian $p$-group.

\begin{definition} If $V$ is a torsion-free divisible group we say we are in characteristic 0. If $V$ is an elementary abelian $p$-group for some prime $p$, then we say that characteristic is $p$.
\end{definition}

\begin{theorem} If characteristic is $0$, then $G=\operatorname{GL}_n(F)$ and $V=F^n$.
\end{theorem}

\proof When $F$ is an \acf of characteristic~0, any definable additive endomorphism of $F^n$ is an $F$-linear map \cite[Corollary 3.3]{poizat}. Therefore, $G$ can be viewed as a subgroup of $\operatorname{GL}(F^n)=\operatorname{GL}_n(F)$ and hence the result follows from Fact~\ref{abclin}. \hfill $\Box$

\begin{proposition}\label{even}
If characteristic is $2$, then $G$ is of even type; that is, Sylow$^\circ$ $2$-subgroups of $G$ are non-trivial unipotent groups.
\end{proposition}

\proof The Sylow$^\circ$ 2-subgroup of $G$ is of the form $T*B$, where $T$ is a 2-torus and $B$ is a unipotent 2-group  \cite[Proposition I.6.4]{abc}. Since we are in characteristic 2, $V$ is an elementary abelian 2-group, hence both $V$ and $T$ are locally finite. Thus $TV$ is locally finite, and hence $TV$ is nilpotent-by-finite \cite[Proposition I.5.28]{abc}. By \cite[Proposition I.5.8]{abc}, $T$ centralizes $V$, therefore $T=1$ by faithfullness. If $B=1$ (we say $G$ is degenerate when $B=T=1$), then the subgroup generated by pseudoreflection subgroups of $G$ is a normal abelian subgroup \cite[Proposition IV.5.2]{abc}. By Fact~\ref{zilberlin} and Fact~\ref{abclin}, $G\cong GL(V)$, since the characteristic of the underlying field is 2, it follows that $B\neq 1$.   \hfill $\Box$

In characteristic 2, we have the result after combining the following fact with the main result of \cite{abc}: all groups of even type are $K$-groups.

\begin{fact}[Theorem III.1.5 \cite{abc}] Let $V$ be an elementary abelian $2$-group, $G$ a connected $K$-group of even type acting on $V$ definably, faithfully, irreducibly and with a pseudoreflection subgroup. Then there exists an \acf $F$ of characteristic $2$, such that $V$ is a vector space over $F$ and\/ $G\cong \operatorname{GL}(V)$, and the action is the natural action of\/ $\operatorname{GL}(V)$ on $V$.
\end{fact}

\medskip
As the results of this section show we may exclude certain cases from our further analysis.

\medskip
\noindent
{\bf Conclusion.} In this work, we may assume $\rk(V)\geqslant 3$, $Z(G)$ is finite and the characteristic is at least 3.

\section{Some Preliminary Results}
\label{prelim}

Let us recall the simple but important corollary of Zilber's Theorem (Fact~\ref{zilberabelianaction}): if $R$ is a pseudoreflection subgroup in $G$, then $R\cong F^*$ and $[V,R]\cong F^+$ for some \acf $F$.

A direct product of pseudoreflection subgroups in $G$ is obviously a torus; we call it a {\em pseudoreflection torus}. Thus a pseudoreflection torus is isomorphic to a direct product of finitely many copies of $F^*$, hence it is a good torus. We define the {\em pseudorank} of $G$ to be the number of copies of $F^*$ in a maximal pseudoreflection torus of $G$, we denote it by $\psrk(G)$.

We work under the following assumptions:  $G$ is a connected group of \fmr acting on a
connected abelian group $V$ of \fmr definably, faithfully and
irreducibly. Also $G$ contains a pseudoreflection subgroup $R$,
that is a connected definable abelian subgroup such that $V=[V,R]\oplus
C_V(R)$,  and $R$ acts transitively on the nonzero elements of
$[V,R]$. Moreover, we assume $\rk[V,R]=1$, $\rk(V)=n\geqslant 3$, $Z(G)$ is finite and the characteristic is at least 3.

\begin{proposition}\label{odd}
If characteristic is not $2$, then $G$ is of odd type; that is, Sylow$^\circ$ $2$-subgroups of\/ $G$ are $2$-tori.
\end{proposition}

\proof Let the Sylow$^\circ$ 2-subgroup of $G$ be $T*B$, where $T$ is a 2-torus and $B$ is a unipotent 2-group. Since we are not in characteristic 2,  being a unipotent 2-group, $B$ acts trivially on $V$ \cite[Proposition I.8.5]{abc}. Hence by faithfullness $B=1$. Also, note that $R\cong F^*$ has an infinite 2-subgroup. Hence, $T\neq 1$, and $G$ is of odd type. \hfill $\Box$

\begin{lemma} Let $D$ be a pseudoreflection torus in $G$, say $D=R_1\times\cdots\times R_m$ for some pseudoreflection subgroups $R_i$ in $G$. Then $$[V,D]=\bigoplus_{i=1}^m [V,R_i].$$

\end{lemma}

\proof Let $R$ and $S$ be two commuting pseudoreflection subgroups in $G$. It suffices to show that $[V,R]\cap[V,S]\neq 0$ implies $R\cap S\neq 1$.

Assume there is a non-zero element $w\in[V,R]\cap[V,S]$, then we can write $$w=\sum_{i=1}^n(r_iv_i-v_i)=\sum_{i=1}^m(s_iu_i-u_i)$$ for some $r_i\in R$, $s_i\in S$, $v_i,u_i\in V$.

To show $[V,R]=[V,S]$, take a generator $px-x\in [V,R]$, then by transitivity, there exists $t\in R$ such that
\begin{eqnarray*}
px-x &=& t\left(\sum_{i=1}^n(r_iv_i-v_i)\right)\\
    &=& t\left(\sum_{i=1}^m(s_iu_i-u_i)\right)\\
    &=& \sum_{i=1}^m(s_itu_i-tu_i)\\
    &=& \sum_{i=1}^m(s_iw_i-w_i)
\end{eqnarray*}
after setting $w_i=tu_i$. Hence, $px-x\in[V,S]$, that is $[V,R]\subseteq [V,S]$. By symmetry, we get the equality.

To show $C_V(R)=C_V(S)$, let $i_R$ and $i_S$ denote the unique involutions in $R$ and $S$, respectively. Take an arbitrary $v\in C_V(R)$, then one can write $v=u+w$, where $u\in [V,S]=[V,R]$ and $w\in C_V(S)$. Now $i_Si_R(v)=i_S(v)=-u+w$, and on the other hand $i_Si_R(v)=i_Ri_S(v)=i_R(-u+w)$. Thus, $-u+w\in C_V(i_R)=C_V(R)$, hence $u\in [V,R]\cap C_V(R)=0$, that is $v=w\in C_V(S)$. By symmetry, the result follows.

Above claims show that the actions of $i_R$ and $i_S$ are identical on $V$,
thus by faithfullness $i_R=i_S$ and hence $R\cap S\neq 1$.%

\hfill $\Box$

\begin{lemma} Every $2$-element in $G$  {\rm (}and also in $G/Z${\rm )} belongs to a $2$-torus.\label{2toral}
\end{lemma}

\proof  Apply Proposition~\ref{odd} and Fact~\ref{burche} with $p=2$. \hfill $\Box$

\begin{lemma} For every pseudoreflection subgroup $R\leqslant G$, $C_G(R)=C_G^\circ(R)=R\times K_R$, where $K_R$ is the kernel of the action of $C_G(R)$ on $[V,R]$. Moreover, $K_R$ acts faithfully on $C_V(R)$. More generally, for a pseudoreflection torus $D=R_1\times\cdots\times R_m$,
$C_G(D)=D\times K_{m}$, where $K_{m}$ acts faithfully on $\bigcap_{i=1}^m C_V(R_i)$. \label{centroftori}

\end{lemma}

\proof Since $R\cong F^*$, it is a good torus, hence its centraliser is connected. It is clear that $C_G(R)$ normalizes $[V,R]$. Let $K_R$ be the kernel. Obviously, $R\times K_R\leqslant C_G(R)$. Then $C_G(R)/K_R$ acts faithfully and commuting with the action of $R\cong F^*$ on $[V,R]$. Hence, $C_G(R)/K_R$ acts linearly on $[V,R]\cong F^+$, that is, it is one-dimensional. Hence the first two statements are proved. The moreover part is easy and the last statement can be proven similarly. \hfill $\Box$

\medskip
From now on, we assume $\psrk(G)=\rk(V)=n\geqslant 3$. Let $T=R_1\times\cdots\times R_n$ be a maximal pseudoreflection torus in $G$, where each $R_i$ is a pseudoreflection subgroup. %In this case, note that
%$$\bigcap_{i=1}^j C_V(R_i)=\bigoplus_{i=j+1}^n [V,R_i].$$
%To simplify notation, set $$\hat{R}_{ij}=R_1\times\cdots\times \hat{R_i}\times\cdots\times\hat{R_j}\times\cdots\times R_n.$$

%\begin{corollary}
%For every $i\neq j$, $$C_G(\hat{R}_{ij})=\hat{R}_{ij}\times K_{i,j},$$ where $K_{i,j}$ acts on $[V,R_i]\oplus [V,R_j]$ faithfully, and
%$K_{i,j}=\operatorname{GL}_2$ or is solvable. \label{rk2}
%\end{corollary}

%\proof Use Fact~\ref{deloro} and note that $\psrk(K_{i,j})=2$. \hfill $\Box$

At this point, we would like to emphasize that maximal 2-tori were central in the study of groups of odd type. Hence, we make the following observation on 2-tori in $G$.

\begin{lemma} Let $T$ be a maximal pseudoreflection torus in $G$. Then
$C_G(T)=C_G^\circ(T)=N_G^\circ(T)=T$, and hence $T$ is a maximal torus in $G$. In particular, every maximal pseudoreflection torus in $G$ contains a maximal $2$-torus of $G$, and $Z(G)$ lies in every maximal torus of $G$. \label{torus}
\end{lemma}

\proof Since $T\cong\bigoplus F^*$, $T$ is a good torus. Thus, the first two equalities follow from Fact~\ref{cherlintorus}. Note $R_1\times\cdots\times R_{n-1}\leqslant T$, thus $T\leqslant C_G(T)\leqslant C_G(R_1\times\cdots\times R_{n-1})=R_1\times\cdots\times R_{n-1}\times R_n=T$ by Lemma~\ref{centroftori}.
Hence we are done. \hfill $\Box$

\begin{lemma} The center $Z(G)$ is a non-trivial finite cyclic group
and contains a unique involution.
\label{center}
\end{lemma}

\proof Recall that we assume $Z=Z(G)$ is finite.

Since $G$ acts irreducibly on $V$, the ring $R$ of definable endomorphisms of $V$ is a divison ring by Schur's Lemma (\cite[Lemma I.4.7]{abc}). Hence, the subring $S$ generated by $Z$ in $R$ is a definable integral domain (\cite[Proof of Lemma I.4.8]{abc}), and is contained in a definable division ring, which is a field \cite[I.4.25]{abc}. Being a finite subgroup of the multiplicative subgroup of a field, $Z$ is cyclic.

By Lemma~\ref{torus}, $Z$ lies in every torus, and in particular in $T=R_1\times\cdots\times R_n$, so we will represent elements of $Z$ as $(r_1,\ldots,r_n)$. Since $(-1,\ldots,-1)$ inverts all the elements in $V$, clearly $Z\neq 1$. \hfill $\Box$

\begin{lemma}\label{G/Z} The group $G/Z$ is centerless, and has more than one conjugacy classes of involutions.  In particular, $G/Z$ has no strongly embedded subgroups.
\end{lemma}

\proof By Lemma I.3.8 \cite{abc}, $G/Z$ has no center.
Let the order of $Z(G)$ be $2^ml$ where $l$ is odd. Since $Z(G)$ has a unique involution, $Z(G)$
contains a cyclic subgroup of order $2^m$. Let $d$ be the generator of this cyclic subgroup. Since $d\in Z(G)$ lies in every maximal torus, by divisibility of tori, for every maximal torus $T$, there exists $d_T\in T$ such that $d_T^2=d$.
Then $d_T$ is a non-central 2-element of order $2^{m+1}$ where $\bar{d_T}$ is an involution.
Let $j$ be the unique involution in $R_1$.
Assume $\bar{d_T}$ and  $\bar{j}$ are conjugate in $\bar{G}$, say
$h^{-1}jh=d_Tz$ for some $h\in G$ and $z\in Z$. Then $1=(h^{-1}jh)^{2^m}=d_T^{2^m}z^{2^m}=(-1)z^{2^m}$, and hence $z^{2^m}=-1$. Thus $z$ is an element of order $2^{m+1}$ in $Z(G)$, which is a contradiction. Therefore, $\bar{d_T}$ and $\bar{j}$ belong to two distinct conjugacy classes of involutions in $G/Z$. The last statement follows from \cite[I.10.12]{abc}. \hfill $\Box$

\section{Reductivity}

\label{reduct}

{\bf New Assumption.} From now on we will work with a minimal counterexample to Theorem~\ref{main}: we assume Theorem~\ref{main} is false, and among all counterexamples first we pick one with minimum $\rk(V)$; then with minimum $\rk(G)$. Thus, $\rk(V)\geqslant 3$ and $Z(G)$ is finite. Let us denote the exponent of $V$ by $p$ and note that we can assume $p\geqslant 3$.

\begin{lemma} Let\/ $H$ be a connected definable subgroup of\/ $G$. Then $H$ is unipotent {\rm (}that is, nilpotent of bounded exponent{\rm )} if and only if\/ $H \ltimes V$ is nilpotent.
\end{lemma}

\proof Since $V$ is abelian and $H$ is nilpotent, $HV$ is solvable. By \cite[I.8.4]{abc}, $HV\leqslant F^\circ(HV)$, thus $HV$ is nilpotent. Conversely, if $HV$ is nilpotent, we can write $H=UT$, where $U$ is unipotent and $T$ is radicable. Since $HV$ is nilpotent, $HV=T(UV)$ is its decomposition, since $UV$ is of bounded exponent. By  \cite[I.8.4]{abc}, $T$ commutes with $UV$, hence $T=1$ by faithfulness. So $H=U$ is unipotent. \hfill $\Box$

\begin{proposition} Let $H$ be a definable subgroup of\/ $G$. Assume that
$$0=V_0 < V_1 < \cdots < V_l = V$$
is a composition series for the action of $H$ on $V$. Then the kernel
of the natural action of\/ $H$ on
$V_1/V_0 \oplus \cdots \oplus V_l/V_{l-1}$
is a unipotent $p$-group.
\end{proposition}

\proof The kernel $K$ of the above action is a nilpotent subgroup by \cite[16.3.1]{kargapolov}. Say $K=TU$, where $T$ is the radicable part and $U$ is unipotent \cite[I.5.8]{abc}. Since $T$ contains no infinite $p$-unipotent subgroup, by \cite[I.8.5]{abc}, $T$ centralizes $V$, and hence $T=1$ by faithfulness. Thus, $K=U$ is unipotent. The fact that $K$ is a $p$-group follows from \cite[I.5.16]{abc}. \hfill $\Box$

\definition Denote this kernel by $R_u(H)$ and call it the
\emph{unipotent radical} of $H$.
In fact, it follows from the discussion in the previous paragraph that $R_u(H)$ is the maximal definable
normal unipotent subgroup of $H$.

\begin{lemma} Let $H$ be a proper connected definable subgroup in $G$ containing $T$. Then
$$H/R_u(H) \cong GL_{n_1}(F) \oplus\cdots\oplus GL_{n_k}(F)$$
for some $n_1+\cdots +n_k = {\rm rk}(V)$. \label{reductive} \end{lemma}

\proof Recall that $\rk(V)\geqslant 3$. First note that if $H$ acts irreducibly on $V$, then by induction assumption $H\cong GL(V)$, hence we are done. So, let's assume that $H$ acts reducibly on $V$ with a composition series $0=V_0<V_1<\cdots< V_k=V$. Let $H_i$ be the kernel of the action of $H$ on $V_i/V_{i-1}$, for $i=1,\dots, k$. Then $\bigcap H_i=R_u(H)$ and hence $H/R_u(H)$ embeds into $H/H_1\times\cdots\times H/H_k$.
Note that $H/H_i$ acts faithfully and irreducibly on $V_i/V_{i-1}$. Since for every $i$, $H/H_i$ contains a pseudoreflection subgroup, by induction assumption, each $V_i/V_{i-1}\cong F^{n_i}$ and $H/H_i\cong GL_{n_i}(F)$. Therefore, $V=F^{\sum n_i}$ and $H/R_u(H)\leqslant \bigoplus GL_{n_i}$. Now note $H/R_u(H) \cong \bigoplus GL_{n_i}$.
Each $R_i\leqslant T$ maps into exactly one of the direct summands (isomorphic to $GL_{n_i}$), hence the diagonal subgroup of $\bigoplus GL_{n_i}$ lies in the image of the embedding. Hence, each $H/H_i$ embeds into one direct summand, thus the embedding is onto. Note that $\rk(V)=(\sum n_i)\rk (F)=(\sum n_i)\rk[V,R]=\sum n_i$. \hfill $\Box$

\medskip

Set $Z=Z(G)$, since $Z$ is finite $G/Z$ is a centerless group (Lemma~\ref{G/Z}). Use $\bar{H}=H/Z$, if $Z\leqslant H\leqslant G$. Note that $Z$ lies in every maximal torus, by Lemma~\ref{torus}. Also set $R_u(\bar{H})=R_u(H)Z/Z\cong R_u(H)$, which is nilpotent.

\begin{corollary}\label{quotientreductive} Let $\bar{H}$ be a proper connected definable subgroup in $\bar{G}$ containing $\bar{T}$. Then
$$\bar{H}/R_u(\bar{H}) \cong (GL_{n_1}(F) \oplus\cdots\oplus GL_{n_k}(F))/A,$$ where $A$ is a finite central subgroup and
$n_1+\cdots +n_k = {\rm rk}(V)$.
\end{corollary}

\proof This follows immediately from Lemma~\ref{reductive}. Consider $\bar{H}/R_u(\bar{H})=(H/Z)/(R_u(H)Z/Z)\cong H/R_u(H)Z\cong (H/R_u(H))/(R_u(H)Z/R_u(H))$. Note that by Lemma~\ref{reductive}, $H/R_u(H)\cong\bigoplus GL_{n_i}$, and $$R_u(H)Z/R_u(H)\cong Z/(R_u(H)\cap Z)\cong Z$$ is a finite central group. Hence, $\bar{H}/R_u(\bar{H})$ is a quotient of $\oplus GL_{n_i}$ by a finite central subgroup. \hfill $\Box$

\begin{lemma}\label{simple} $G/Z$ is simple.
\end{lemma}

\proof Let $S$ be a minimal normal definable subgroup in $G$. If $S$ is abelian then we are done by Fact~\ref{zilberlin}. Therefore, $S$ is quasisimple \cite[Chapter 7]{bn}. If $ST$ is a proper subgroup in $G$, then by Lemma~\ref{reductive}, $ST/R_u(TS)$ is isomorphic to a direct sum $D$ of $GL_{k_i}$'s where $k_1+\cdots+k_r=n$. Set $R=R_u(TS)$. Then $SR/R\cong S/(S\cap R)$ is a non-abelian quasisimple normal subgroup in $D$, hence $SR/R\cong SL_{k}$ for some $k\geqslant 2$. Therefore,
$ST/R\cong (GL_k)\oplus T_{n-k}$ where $T_i$ stands for an $i$-dimensional torus. Since $S$ is normal in $G$, $[S,V]=V$ and $C_V(S)=0$, thus $n-k=0$ that is $ST/R\cong GL(V)$. In particular, $ST$ acts irreducibly on $V$, by induction hypothesis, $ST\cong GL(V)$, $V=F^n$ and $ST$ acts linearly on $V$.

Now set $N=\langle R^g\mid g\in G\rangle$. Then  $ST\leqslant N\unlhd G$. Now note that $ST=N$. Indeed, an arbitrary $R^g\leqslant N$ acts trivially on $g^{-1}C_V(R)$ and acts like $F^*$ on $g^{-1}[V,R]\cong F^+$. However, given a direct decomposition of the vector space $F^n$ into two subspaces $A$ and $B$, of dimension 1 and $n-1$, respectively, there always exists a subgroup $D\cong F^*$ in $GL_n(F)$ that centralizes $B$ and acts like $F^*$ on $A$. By the faithfullness of the action, $R^g\leqslant ST$, and hence $N=ST$.
Thus the infinite center of $N=GL(V)$ is normal in $G$. Again, we are done by Fact~\ref{zilberlin}.

Now assume $G=ST$. Note  $T=(T\cap S)^\circ\oplus T_0$ for some subtorus $T_0\leqslant T$, therefore $T_0$ centralizes $S$, and hence $G$, and we get $T_0=1$, that is, $T\leqslant S$ and hence $G=S$ is quasisimple.  \hfill $\Box$

\begin{theorem} For every involution $\bar{i}\in \bar{G}$, $$C_{\bar{G}}^\circ(\bar{i})\cong (GL_{n_1}(F) \oplus\cdots\oplus GL_{n_k}(F))/A,$$ where $A$ is a finite central subgroup and $n_1+\cdots +n_k = {\rm rk}(V)$. \label{centrs}
\end{theorem}

\proof  Let $\bar{i}$ be an involution in $\bar{G}$, then by Lemma~\ref{2toral}, $\bar{i}$ lies in a torus, say $\bar{T}$.
First note that $\bar{i}\in \bar{T}\leqslant C^\circ_{\bar{G}}(\bar{i})$. If we assume $R_u(C^\circ_{\bar{G}}(\bar{i}))\neq 1$, then we can define a non-trivial nilpotent signalizer functor on the set of involutions in $\bar{G}$ by setting $$\theta(\bar{j})= R_u(C^\circ_{\bar{G}}(\bar{j})).$$
Indeed by Corollary~\ref{quotientreductive},
this functor satisfies the Balance Property, that is $$\theta(\bar{j})\cap C_{\bar{G}}(\bar{i})\leqslant \theta(\bar{i}).$$

We are now in position to use the Uniqueness Subgroup Theorem \cite{bbn2}.

\begin{fact}[Uniqueness Subgroup Theorem]
Let $G$ be a simple group of finite Morley rank and odd type with
Pr\"{u}fer $2$-rank $\geqslant 3$.  Let $S$ be a Sylow
$2$-subgroup of\/ $G$ and $D=S^\circ$ the maximal\/ $2$-torus in $S$. Let further $E$ be the subgroup of\/ $D$ generated by involutions and $\theta$ a nilpotent signaliser functor on $G$.

Assume in addition that every proper definable subgroup of\/ $G$ containing $T$ is a $K$-group.

Then $M = N_G(\theta(E))$ is a  strongly embedded subgroup in $G$.

\label{f:SET}
\end{fact}

Thus, $\bar{G}$ has a strongly embedded subgroup, and this contradicts with Lemma~\ref{G/Z}. Therefore,  $R_u(C^\circ_{\bar{G}}(\bar{i}))=1$, and by Corollary~\ref{quotientreductive}, the structure of centralisers are as described above. \hfill $\Box$

\section{The Generic Identification Theorem}
\label{sec;git}

In this section, we will show that $G/Z$ satisfies the conditions of the Generic Identification Theorem \cite{git,git2} and hence we will conclude that $G/Z$ is a simple Chevalley group over an \acf of characteristic $p\neq2$, therefore $G$ is a quasisimple Chevalley group over an \acf of charactersitic $p\neq2$.

In fact, the original theorem is stated for every prime, but we will use it for the prime 2 only, so we state it in that special form. We need a version of the theorem with slightly weaker assumptions than that of \cite{git}; essentially the same proof as in \cite{git} works, and this is documented in \cite{git2}.

\begin{fact}[Generic Identification Theorem \cite{git2}]
Let\/ $G$\/ be a simple  group of finite
Morley rank and\/ $D$ a maximal\/ $2$-torus  in $G$
of Pr\"{u}fer rank at least\/ $3$. Assume that
\begin{itemize}
\item[{\rm (A)}] Every proper connected definable subgroup of $G$ which contains $D$ is a $K$-group.

\item[{\rm (B)}] For every involution\/ $x$ in $D$, the group $C^\circ_G(x)$ contains no infinite elementary abelian $2$-subgroup and
$$C^\circ_G(x)=F^\circ(C^\circ_G(x))E(C^\circ_G(x)).$$

\item[{\rm (C)}] $\langle C^\circ_G(x) \mid x \in D, \;|x|=2\rangle = G$.
\end{itemize}
Then $G$ is a Chevalley group over an \acf of characteristic distinct
from $2$.

\label{siskebap}
\end{fact}

Next, we will show that $G/Z$ satisfies the conditions of the Generic Identification Theorem. First, note that $G/Z$ is simple by Lemma~\ref{simple}. We will take $D$ to be the Sylow 2-subgroup of $T$, note that the definable hull of $D$ is $T$. Since
all definable subgroups of  $G/Z$  containing a maximal pseudoreflection torus are $K$-groups as discussed before, (A) is satisfied.

From now on we work in $G/Z$, and we write $\bar{S}$ for the image of the natural epimorphsim $G\to G/Z$ for every subset $S$ in $G$.

\begin{itemize}
\item[(B)]
For every involution $\bar{x}\in \bar{D}$, $C^\circ_{\bar{G}}(\bar{x})$ contains no infinite elementary abelian 2-subgroup and $C^\circ_{\bar{G}}(\bar{x})=F(C^\circ_{\bar{G}}(\bar{x}))^\circ E(C^\circ_{\bar{G}}(\bar{x}))$.
\end{itemize}

\proof Since $G$ is of odd type by Lemma~\ref{odd}, so is $\bar{G}$. Hence, no subgroup of $G$, or $\bar{G}$, contains an infinite elementary abelian 2-subgroup. The second part of the statement follows from Theorem~\ref{centrs}.

\begin{itemize}
\item[(C)] $\bar{G}=\langle C^\circ_{\bar{G}}(\bar{x}) \mid \bar{x} \in \bar{D}, \;|\bar{x}|=2\rangle$.
\end{itemize}

\proof Set $$\bar{L}=\langle C^\circ_{\bar{G}}(\bar{x}) \mid \bar{x} \in \bar{D}, \;|\bar{x}|=2\rangle.$$
We will show that $N_{\bar{G}}(\bar{L})$ is a strongly embedded subgroup in $\bar{G}$. Then Lemma~\ref{G/Z} implies $\bar{G}=N_{\bar{G}}(\bar{L})$, and then by the simplicity of $\bar{G}$, $\bar{G}=\bar{L}$ will follow.

Let $\bar{S}$ be a Sylow $2$-subgroup in $\bar{G}$ containing $\bar{D}$.
Observe that $N_{\bar{G}}(\bar{S})\leqslant N_{\bar{G}}(\bar{S}^\circ)\leqslant N_{\bar{G}}(\bar{T})$, since the definable hull of $S^\circ$ is $T$. Let $\bar{a}\in N_{\bar{G}}(\bar{S})$ and $\bar{x}$ be an involution in $\bar{D}$. Since $\bar{x}^{\bar{a}}\in \bar{D}$, we get $C^\circ_{\bar{G}}(\bar{x})^{\bar{a}}$ lies in $\bar{L}$. Therefore, $N_{\bar{G}}(\bar{S})\leqslant N_{\bar{G}}(\bar{L})$.

Now let $\bar{i}$ be an involution in $\bar{S}$. Then $C^\circ_{\bar{G}}(\bar{i})$ is a $K$-group by Theorem~\ref{centrs}. On the other hand, since $\bar{i}$ normalizes $\bar{S^\circ}=\bar{D}$, and $\bar{D}$ has Pr\"ufer rank at least 3; the Pr\"ufer rank of $C_{\bar{D}}(\bar{i})$ is at least 2. Therefore, $C_{\bar{D}}(\bar{i})$ contains a 4-group, say $V$. Then, by \cite[Theorem~5.14]{nato},
$$C^\circ_{\bar{G}}(\bar{i})=\langle C_{C^\circ_{\bar{G}}(\bar{i})}({v})\mid v\in V, v\neq 1\rangle.$$
Therefore, $C^\circ_{\bar{G}}(\bar{i})\leqslant \bar{L}$, since $V\leqslant \bar{D}$.

This shows that $N_{\bar{G}}(\bar{L})$ is a strongly embedded subgroup in $\bar{G}$, and hence (C) is satisfied, as explained above.

\bigskip

Therefore, $G$ is a quasisimple Chevalley group of Lie rank at least $3$ over an \acf of charactersitic $p$ and $p\geqslant 3$. The centralisers of involutions in these groups are well-known, see for example \cite{gls}. In our case, the centralisers of involutions are direct sums of copies of $\operatorname{GL}_{k_i}$. We compare the centralisers in Table~\ref{table}.

\begin{table}[!ht]
\begin{center}
\begin{tabular}{ccc}
\\ % \toprule
Type & $n$ & $C_G(i)$ for some $i\in G$
\\ \\ % \midrule

$A_n$ & $n\geqslant 3$ & $A_2A_{n-2}T_1$\\

$B_n$ & $n\geqslant 3$ &  $B_{n-1}T_1$\\

$C_n$ & $n\geqslant 3$ & $C_1C_2$\\

$D_4$ &   &  $(A_1)^4$\\

$D_n$ & $n\geqslant 5$ & $D_{n-1}T_1$\\

$E_6$ &  & $D_5T_1$\\

$E_7$ &  & $D_6A_1$\\

$E_8$ & &  $D_8$\\

$F_4$ & & $A_1C_3$

\end{tabular}
\end{center}
\caption{Some centralisers of involutions in Chevalley groups over algebraically closed fields of odd or zero characteristic.}
\label{table}
\end{table}

Clearly, none of the centralisers in the table is isomorphic to a direct sum of copies of $GL_{n_i}$.
This final contradiction shows that there is no counterexample to Theorem~\ref{main}. \hfill $\Box$

\section{Large 2-torus}
\label{sec:large-torus}

In this section we shall prove Theorem~\ref{large-torus}.
We start with analysing a faithful definable action of a ``large'' elementary abelian $2$-group $E$ of order $2^m$ on a connected abelian group $V$ (written additively) of finite Morley rank which has odd prime exponent $p$ or is divisible and torsion-free. Assume that $\rk V = n$. The ``largeness'' of $E$ will mean that $m \geqslant n$. As we shall soon see, these assumptions will lead to a very concrete and explicit configuration.

We need to use some elementary standard concepts from representation theory. A \emph{character} of $E$ is a homomorphism $\rho: E \to \{\pm 1\}$; the pointwise multiplication of characters turns the set of characters into a group $E^*$ called the \emph{dual group} of $E$; it is well-known that since $E$ is a finite elementary abelian $2$-group, $E \simeq E^*$.

\begin{lemma}\label{lm:2-group}
Let $V$ be a connected abelian group and $E$ an elementary abelian 2-group of order $2^m$ acting definably and faithfully on $V$. Assume $m\geqslant n=\rk(V)$ and $V$ is either torsion-free divisible, or of exponent an odd prime. Then $m = n$ and $V = V_1 \oplus\cdots\oplus V_n$, where
\bi
\item[{\rm (a)}]
every subgroup $V_i$, $i = 1,\dots,n$, is connected, has Morley rank $1$ and is $E$-invariant.
\ei
Moreover,
\bi
\item[{\rm (b)}] for each $V_i$,  $i = 1,\dots,n$, is a weight space of $E$, that is, there exists a character $\rho_i \in E^*$ such that
     \[
     V_i = \{ v \in V \mid v^e = \rho_i\cdot v \mbox{ for all }   e \in E\}.
     \]
\ei
\end{lemma}

\proof
Observe further that if $e\in E$ is an involution, $V = C_V(e) \oplus [V,e]$ and both subgroups $C_V(e)$ and $[V,e]$ are definable, connected and $E$-invariant. Consider a direct decomposition $V = V_1 \oplus\cdots\oplus V_k$ into a direct sum of non-trivial connected definable $E$-invariant subgroups with the maximal possible number $k$ of summands; observe that for each $i$ we have $\rk(V_i)\geqslant 1$ and therefore $k \leqslant n$. Then, given an involution $e \in E$ and arbitrary $V_i$, $e$ either centralises $V_i$ or acts on $V_i$ by inversion: $v^e = -v$ for all  $v \in V_i$. Therefore for each $V_i$, $i= 1,\dots, k$ there exists a character $\rho_i \in E^*$ such that $v^e = \rho_i\cdot v $ for all $v \in V_i $ and $e \in E$. Since the action of $E$ on $V$ is faithful, the map
\[
e \mapsto (\rho_1(e),\dots,\rho_k(e))
\]
is an embedding of $E$ into $\{\pm 1\}^l$ for $l \leqslant k$ (we can here into account that some of $\rho_i$ and $\rho_j$ could be equal) and therefore $m \leqslant l \leqslant k\leqslant n$; combining that with the inequality $m \geqslant n$ we have $m=l=k=n$. In particular, this means that all $i$ characters $\rho_i$ are distinct, $\rk (V_i) =1$ and $V_i$ are weight spaces for $V$. This completes the proof of the lemma. \hfill $\Box$

\begin{lemma}
Assume that a connected group of finite Morley rank\/ $G$ acts faithfully and definably on a connected abelian group $V$, where $V$ is either torsion-free divisible or of exponent an odd prime. Let\/ $D$ be a maximal $2$-torus in $G$ of Pr\"ufer $2$-rank at least\/ $n=\rk(V)$, then the definable closure of\/ $D$ is a pseudoreflection torus of pseudoreflection rank $n$.
\end{lemma}

\proof We will use induction on $n=\rk(V)\geqslant 1$. If $n=1$, then we can apply Lemma~\ref{strmin}, and get an \acf $F$, where $V\cong F^+$
and the definable closure $T$ of $D$ lies in $F^*$. Since $\rk(F^*)=1$ and the action is the usual multiplication, we conclude that $T\cong F^*$ is a pseudoreflection torus of pseudoreflection rank 1.

Assume $n\geqslant 2$, and let $E$ be the subgroup generated by all involutions in $D$. Then $|E|\geqslant 2^n$ and we are in the configuration described in Lemma~\ref{lm:2-group}; retain the notation of that lemma. Let $T$ stand for the definable closure of $D$, then $T$ acts non-trivially on connected groups $V_i= \{v\in V \mid v^{e_i}=-v\}$. Write $W_i = \bigoplus_{j\ne i} V_j$ and $R_i = C_T(W_i)$ for every $1\leqslant i\leqslant n$. Then for each $1\leqslant i\leqslant n$, $T/R_i$ acts faithfully on $W_i$. Since $W_i$ is of rank $n-1$, the Pr\"ufer rank of $T/R_i$ is at most $n-1$. By induction hypothesis, $T/R_i$ is a pseudoreflection torus of rank $n-1$, thus $R_i$ is an infinite definable group acting on $V_i$ faithfully. By Lemma~\ref{strmin}, $R_i$ is the multiplicative group $F^*$ of the field arising from its action on $V_i$ treated as an additive group of $F$, while centralising $W_i$. Now
\[
T = R_1 \times \cdots \times R_n
\]
is a pseudoreflection torus of pseudoreflection rank $n$. \hfill $\Box$

\bigskip

Theorem~\ref{large-torus} now follows from Theorem~\ref{main}.

\end{document}